\documentclass[12pt]{amsart}
\newcommand{\RP}{{\mathbb {RP}}}
\newcommand{\CP}{{\mathbb {CP}}}
\newcommand{\Z}{{\mathbb Z}}
\newcommand{\R}{{\mathbb R}}

\newcommand{\Q}{{\mathbb Q}}
\newcommand{\sm}{\setminus}

\newtheorem{theorem}{Theorem}

\newtheorem{proposition}{Proposition}
\newtheorem{corollary}{Corollary}
\newtheorem{remark}{Remark}

\newtheorem{lemma}{Lemma}

\begin{document}

\author{V.A.~Vassiliev}
\address{Steklov Mathematical Institute of Russian Academy of Sciences;
National Research Institute --- Higher School of Economics} \email{vva@mi.ras.ru}

\thanks{
Supported by Russian Science Foundation under the grant 14-50-00005}

\title{Rational homology of the order complex of zero sets of homogeneous
quadratic polynomial systems in $\R^3$}

\begin{abstract}
The naturally topologized order complex of proper algebraic subsets in $\RP^2$,
defined by systems of quadratic forms, has rational homology of $S^{13}$.
\end{abstract}

\maketitle

\section{Introduction}

Order complexes of zero sets of polynomial systems of various types are an
important part of the study of (discriminant or resultant) varieties of
non-generic systems, and hence of their complementary sets of generic systems
(or other similar objects). For some examples and applications see
\cite{kotor}, \cite{novikov}, \cite{Gorinov}, \cite{Tommasi}. In this article
we study the (naturally compactified) order complex of proper subsets in
$\RP^2$ defined by systems of homogeneous quadratic polynomials; in particular
we calculate its rational homology group.

\subsection{Main result}

Let $W \simeq \R^6$ be the space of all real quadratic forms in $\R^3$. For any
$i=1, \dots, 5$ consider the Grassmann manifold $G_i(W)$ of $i$-dimensional
subspaces in $W$.

Any collection of points in $\RP^2$ defines a vector subspace in $W$ consisting
of all quadratic forms vanishing on all corresponding lines in $\R^3$. Denote
by $J_i \subset G_{6-i}(W)$ the {\em closure} of the union of all points in
$G_{6-i}(W)$ corresponding to such subspaces of codimension $i$.

The disjoint union $J_1 \sqcup \dots \sqcup J_5$ is a partially ordered set,
with the order relation defined by the incidence of corresponding subspaces.

Consider the join $J_1 * \dots * J_5$, i.e. the union of all simplices of
dimensions $\leq 4$, whose vertices correspond to points of different $J_i$,
and define the {\em order complex} $\Psi \subset J_1 * \dots * J_5$ of our
partially ordered set as the union of those such simplices all whose vertices
are incident to one another. $\Psi$ is a closed subset of the compact space
$J_1 * \dots * J_5$, in particular is also compact.

\begin{theorem}
\label{mthm123} The rational homology group of the order complex $\Psi$ is
equal to $\Q$ in dimensions $0$ and $13$ and is trivial in all other
dimensions.
\end{theorem}

\begin{remark} \rm
Similar order complexes related with different spaces of functions participate
in many known theories and theorems. For instance, by a version of
C.~Caratheodory theorem, the order complex of proper subsets in $\RP^1$
specified by systems of homogeneous polynomials of degree $r$ in $\R^2$ is
homotopy equivalent to $S^{2r-1}$, see \cite{KK}. The order complex of proper
vector subspaces in $\R^n$ (= subsets specified by systems of linear equations)
is homeomorphic to $S^N, N=n(n-1)/2 + n-2$. Also, the homology groups of
certain such complexes related with complex homogeneous polynomials contribute
to the homology groups of moduli spaces of complex curves, see  \cite{Tommasi},
\cite{novikov}, \cite{Gorinov}.
\end{remark}

\begin{remark} \rm
The order complex $\Psi$ is closely related to the geometry of the resultant
varieties in the spaces of quadratic polynomial systems in $\R^3$, in
particular it is the principal part of the construction of one of possible
simplicial resolutions of these varieties. These spaces were studied in
\cite{ser}, where however a different construction of resolutions (the one
proposed in \cite{Gorinov}) was used.
\end{remark}

\noindent {\bf Problems.} There are obvious questions concerning the extension
of our calculation to integral homology groups, and also to polynomial systems
of higher degrees and/or in spaces of higher dimensions.

\subsection{Main spectral sequence}

For any $p=1, \dots, 5$ and any point $X \subset J_p$ define its {\em
subordinate subcomplex} $\Psi(X)$ as the union of those simplices in $J_1
* \dots * J_p$, all whose vertices correspond to subspaces in $W$ containing
the $(6-p)$-dimensional subspace $\{X\}$. By the construction, any $\Psi(X)$ is
a compact contractible topological space: it is a cone with vertex $X$.

Define the filtration
\begin{equation}
\Psi_1 \subset \dots \subset \Psi_5 \equiv \Psi \label{fltr1}
\end{equation}
in the order complex $\Psi$: its term $\Psi_p$ consists of all simplices, all
whose vertices belong to $J_1, J_2, \dots, J_p$ only; in particular $\Psi_1 =
J_1$. In other words, $\Psi_p$ is the union of subcomplexes $\Psi(X)$ over all
$X \in J_1, \dots, J_p$.

Consider the spectral sequence $E_{p,q}^r,$ calculating the group $H_*(\Psi,
\Q)$ and generated by this filtration. Its term $E_{p,q}^1$ is canonically
isomorphic to the group $H_{p+q}(\Psi_p, \Psi_{p-1}; \Q) \equiv
\overline{H}_{p+q}(\Psi_p \sm \Psi_{p-1}, \Q);$ here $\overline{H}_*$ denotes
the {\it Borel--Moore homology group}, i.e. the homology group of the complex
of locally finite singular chains.

\begin{theorem}
\label{prop21} All non-trivial groups $E^1_{p,q}$ of our spectral sequence are
shown in Fig. \ref{mss}. The differential $d^1: E^1_{4,5} \to E^1_{3,5}$ is an
isomorphism.
\end{theorem}

\unitlength=0.80mm \special{em:linewidth 0.4pt} \linethickness{0.4pt}
\begin{figure}
\begin{picture}(65.00,70.00)
\thicklines \put(10.00,10.00){\vector(1,0){50.00}}
\put(10.00,10.00){\vector(0,1){60.00}} \thinlines
\put(18.00,10.00){\line(0,1){55.00}} \put(26.00,10.00){\line(0,1){55.00}}
\put(34.00,10.00){\line(0,1){55.00}} \put(42.00,10.00){\line(0,1){55.00}}
\put(50.00,10.00){\line(0,1){55.00}} \put(5.00,15.00){\makebox(0,0)[cc]{$-1$}}
\put(6.00,21.00){\makebox(0,0)[cc]{$0$}}
\put(14.00,15.00){\makebox(0,0)[cc]{$\Q$}}
\put(5.00,45.00){\makebox(0,0)[cc]{\small $5$}}
\put(30.00,45.00){\makebox(0,0)[cc]{$\Q$}}
\put(38.00,45.00){\makebox(0,0)[cc]{$\Q$}}
\put(5.00,60.00){\makebox(0,0)[cc]{\small $8$}}
\put(46.00,60.00){\makebox(0,0)[cc]{$\Q$}}
\put(14.00,5.00){\makebox(0,0)[cc]{1}} \put(22.00,5.00){\makebox(0,0)[cc]{2}}
\put(30.00,5.00){\makebox(0,0)[cc]{3}} \put(38.00,5.00){\makebox(0,0)[cc]{4}}
\put(46.00,5.00){\makebox(0,0)[cc]{5}}
\end{picture}
\caption{$E^1$} \label{mss}
\end{figure}

Theorem \ref{mthm123} follows immediately from this one. The proof of Theorem
\ref{prop21} takes the rest of the article.

\subsection{Stratification of spaces $\Psi_p \setminus \Psi_{p-1}$ and $J_p$}

There is a natural map $\Psi_p \setminus \Psi_{p-1} \to J_p$: it associates
with any point $\varphi \in \Psi_p \setminus \Psi_{p-1}$ the unique point $X
\in J_p$ such that $\varphi$ belongs to the cone $\Psi(X)$. Therefore any space
$\Psi_p \setminus \Psi_{p-1}$, $p=1, \dots, 5,$ splits into several pieces,
corresponding to the strata of the natural stratification of the space $J_p$.
The latter stratification is elementary: here is the list of its strata.

\begin{lemma} \label{lem3}
The varieties $J_i$, $i=1, \dots, 5$, consist of the smooth strata, whose
points correspond to the subspaces in $W$, formed respectively by polynomials

\noindent $(J_1)$ \quad Vanishing at some point in $\RP^2$.
\smallskip

\noindent $(J_2)$ \quad $(a)$ Vanishing at some two different points in
$\RP^2$.

\qquad $(b)$ Vanishing at some point in $\RP^2$ and having zero derivative
along a fixed direction at this point.
\smallskip

\noindent $(J_3)$ \quad $(a)$ Vanishing at some three points in general
position $($i.e. not in the same line$)$.

\qquad $(b)$ Vanishing on some line in $\RP^2$.

\qquad $(c)$ Vanishing at some two points in $\RP^2$ and having zero derivative
along a fixed direction at one of these points; this direction is different
from the line between these two points.

\qquad $(d)$ Vanishing with all partial derivatives at a point in $\RP^2$.
\smallskip

\noindent $(J_4)$\quad $(a)$ Vanishing at some four points in general position
$($none three in the same line$)$.

\qquad $(b)$ Vanishing on a line and a point outside it.

\qquad $(c)$ Vanishing at three points in general position and having zero
derivative along a fixed direction at one of them; this direction should not
coincide with the line connecting this point with either of other two.

\qquad $(d)$ Vanishing at two points and having zero derivative along fixed
directions at any of these points; these directions should be different from
the line connecting these two points.

\qquad $(e)$ Vanishing on a line and having a singularity at one fixed point in
it.
\smallskip

\noindent $(J_5)$ \quad $(a)$ Vanishing on a smooth non-empty conic.

\qquad $(b)$ Vanishing on two different lines.

\qquad $(c)$ Vanishing on a line with multiplicity 2. \hfill $\Box$
\end{lemma}

\noindent {\bf Notation}. In what follows we denote by $J_p(\alpha)$ the
stratum of $J_p$ listed in the corresponding item $(J_p)(\alpha)$ of this
lemma. Also, we denote by $\Psi_p(\alpha)$ the preimage of $J_p(\alpha)$ under
the natural map $\Psi_p \setminus \Psi_{p-1} \to J_p$. Over any such stratum
$J_p(\alpha)$, this map is a locally trivial fiber bundle
\begin{equation}
\label{ltfb} \Psi_p(\alpha) \to J_p(\alpha).
\end{equation}
The fiber of this bundle over any point $X \in J_p$ is the order complex
$\Psi(X),$ from which its {\em link} $\partial \Psi(X)$ (i.e. the union of all
simplices not containing the vertex point $X$) is removed (the points of this
link correspond to the points of the lower term $\Psi_{p-1}$ of our
filtration).

This difference $\Psi(X) \setminus \partial \Psi(X)$ is denoted by $\breve
\Psi(X)$.
\medskip

By the exact sequence of the pair $(\Psi(X), \partial \Psi(X))$ we have
\begin{equation} \label{esbh}
\overline{H_i}(\breve \Psi(X)) \equiv \tilde H_{i-1}(\partial
\Psi(X)),\end{equation} where $\tilde H_*$ is the homology group reduced modulo
a point.
\medskip

Before studying these spaces $\Psi_p(\alpha)$, we need to calculate or remind
the homology groups of several configuration spaces.

\section{Homology groups of some configuration spaces}

\noindent {\bf Notation.} Given a topological space $M$, denote by $B(M,p)$ the
{\em configuration space}, whose points are the subsets of cardinality $p$ in
$M$. Denote by $\widehat B(\RP^2,p)$ the space of {\em generic}
$p$-configurations in $\RP^2$, i.e. of those subsets of cardinality $p$, none
three points of which belong to one and the same line in $\RP^2$.

For any manifold $M$, we denote by ${{\mathbb O}r}$ the {\it rational
orientation local system} of groups on $M$: it is locally isomorphic to the
constant $\Q$-system, but the elements of $\pi_1(M)$ violating the orientation
act on it as the multiplication by $-1$. For any topological space $M$ denote
by $\pm \Z$ the local system on $B(M,p)$, which is locally isomorphic to the
constant $\Z$-bundle, but the elements of $\pi_1(B(M,p))$ defining odd
permutations of $p$ points act on it as multiplication by $-1$. Also, set $\pm
\Q \equiv \pm \Z \otimes \Q$.

\begin{lemma}[see e.g. \cite{ser}] \label{lem1} For any $p,$ there is a locally
trivial fiber bundle $B(S^1,p) \to S^1$, whose fiber is homeomorphic to
$\R^{p-1}$. This fiber bundle is orientable $($and hence trivial$)$ if $p$ is
odd and is non-orientable if $p$ is even. \hfill $\Box$ \end{lemma}

\begin{lemma}[see e.g. \cite{V}]
\label{lem77} A closed loop in the manifold $B(\RP^2,p)$ violates its
orientation if and only if the cycle in $\RP^2$ composed by the traces of all
$p$ points of our configurations during their corresponding movement defines
the non-trivial element of $H_1(\RP^2, \Z_2)$. \hfill $\Box$
\end{lemma}

\begin{lemma}[see \cite{ser}]
\label{lem11} 1. The Borel--Moore homology groups
$\overline{H}_*(B(\RP^2,2),\Q)$ and $\overline{H}_*(B(\RP^2,2),\pm \Q)$ are
trivial in all dimensions.

2. The group $\overline{H}_i(B(\RP^2,2), {\mathbb O}r)$ is equal to $\Q$ for
$i=4$ and $i=1$ and is trivial in all other dimensions. \hfill $\Box$
\end{lemma}

\begin{lemma}
\label{lem12} The spaces $\widehat B(\RP^2,3)$ and $\widehat B(\RP^2,4)$ are
homotopy equivalent to one another and to the space of unit cubes centered at
the origin in $\R^3$.
\end{lemma}

{\em Proof.} Define the function on $\widehat B(\RP^2,p)$ equal to the minimal
angle between $p$ lines corresponding to the points of configurations. It is
easy to construct smooth vector fields on $\widehat B(\RP^2,3)$ and $\widehat
B(\RP^2,4)$, which strictly increase this function everywhere except for the
manifolds of their maximal points, and thus provide deformation retractions of
these spaces to these manifolds. In the case of $\widehat B(\RP^2,3)$ such a
manifold consists of triples of pairwise orthogonal lines, which can be
considered as normals to the faces of a cube in $\R^3$; so they fix a
homeomorphism between this manifold and the space of all unit cubes. In the
case of $\widehat B(\RP^2,4)$ any similar maximal configuration of four lines
coincides with the set of main diagonals of a cube, which also fixes such a
homeomorphism. \hfill $\Box$

\begin{corollary}
\label{lem17} 1. The manifold $\widehat B(\RP^2,4)$ is orientable, and
$\widehat B(\RP^2,3)$ is not.

2. The orientation sheaf ${\mathbb O}r$ on $\widehat B(\RP^2,3)$ coincides with
$\pm \Q$.

3. The sheaves $\pm \Z$ $($and hence also $\pm \Q)$ on $\widehat B(\RP^2,3)$
and $\widehat B(\RP^2,4)$ are induced from one another by the homotopy
equivalence from Lemma \ref{lem12}.
\end{corollary}

{\it Proof.} By Lemma \ref{lem12}, the fundamental groups of both these
manifolds are generated by $\pi/2$-rotations of a cube around the axes
orthogonal to its faces. Therefore it is enough to check the assertions of our
corollary for any such generator of the fundamental groups. \hfill $\Box$

\begin{lemma}
\label{lem13} 1. The group $\overline{H}_*(\widehat B(\RP^2,3), \Q)$ is trivial
in all dimensions.

2. The group $\overline{H}_*(\widehat B(\RP^2,3), {\mathbb O}r)$ is equal to
$\Q$ in dimensions 6 and 3 and is trivial in all other dimensions.
\end{lemma}

{\em Proof.} 1. Consider the 24-fold covering over $\widehat B(\RP^2,3)$, whose
fiber over a triple of lines consists of {\em positively oriented} ordered
frames of unit vectors in $\R^3$, each of which belongs to one of these lines.
The group $\mbox{GL}_+(3,\R)$ is the direct product of the space of this
covering and $\R_+^3$. Therefore this space is an orientable 6-dimensional
manifold, homotopy equivalent to $\mbox{SO}(3,\R)$; by Poincar\'e duality its
rational Borel--Moore homology groups are equal to $\Q$ in dimensions 6 and 3
and are trivial in all other dimensions. Our covering is finite, therefore (by
Proposition 2 of the Appendix to \cite{Serre}) the homology group
$\overline{H}_*(\widehat B(\RP^2,3), \Q)$ of its base in any dimension is not
greater than that of the total space. So, it is at most 1-dimensional in
dimensions 3 and 6 and is trivial in other dimensions. But
$\overline{H}_6(\widehat B(\RP^2,3), \Q) = 0$ because our space $\widehat
B(\RP^2,3)$ is non-orientable. Hence also $\overline{H}_3(\widehat B(\RP^2,3),
\Q) = 0$ because the Euler characteristic of $\widehat B(\RP^2,3)$ is equal to
0.

2. Also by Poincar\'e duality and Corollary \ref{lem17},
$\overline{H}_*(\widehat B(\RP^2,3), {\mathbb O}r) \simeq H_{6-*}(\widehat
B(\RP^2,3), \Q).$ By the same covering, the latter group can be at most
1-dimensional in dimensions 0 and 3, so the former one is at most 1-dimensional
in dimensions 6 and 3 and is trivial in other dimensions. In dimension 6 a
non-trivial homology class is realized by the fundamental cycle, and the group
$\overline{H}_3(\widehat B(\RP^2,3), {\mathbb O}r)$ should also be non-trivial
by the Euler characteristic considerations. \hfill $\Box$

\begin{corollary} \label{cor22} The manifold of all unit cubes with center $0
\in \R^3$ is orientable.
\end{corollary}

Indeed, by the previous proof and Lemma \ref{lem12} the 3-dimensional rational
homology group of this closed 3-dimensional manifold is non-trivial. \hfill
$\Box$

\begin{lemma}
\label{lem23} 1. The group $\overline{H}_*(\widehat B(\RP^2,4), \Q)$ is equal
to $\Q$ in dimensions 8 and 5 and is trivial in all other dimensions.

2. The group $\overline{H}_*(\widehat B(\RP^2,4), \pm \Q)$ is trivial in all
dimensions.
\end{lemma}

{\it Proof.} 1. By Corollary \ref{lem17} and Poincar\'e duality,
$\overline{H}_*(\widehat B(\RP^2,4), {\mathbb O}r) \simeq H_{8-*}(\widehat
B(\RP^2,4), \Q).$ By Lemma \ref{lem77} and Corollary \ref{cor22}, the latter
group is equal to $\Q$ in dimensions 0 and 3 and trivial in other dimensions.

2. Let $p:I \to \widehat B(\RP^2,4)$ be the 24-fold covering, induced from the
similar covering over $\widehat B(\RP^2,3)$, described in the proof of Lemma
\ref{lem13}, by the homotopy equivalence from Lemma \ref{lem12}. The space $I$
of this covering is orientable, because the base is. By Poincar\'e duality the
group $\overline{H}_i(I, \Q)$ is equal to $\Q$ in dimensions 8 and 5 and is
trivial in other dimensions. The local system $\pm \Q$ is a direct summand of
the local system $p_!(\Q)$ on $\widehat B(\RP^2,4)$. Therefore the homology
group $\overline{H}_i(\widehat B(\RP^2,4), \pm \Q)$ is never greater than
$\overline{H}_i(I,\Q)$, i.e. it can be at most 1-dimensional in dimensions 8
and 5, and is trivial in remaining dimensions. However, $\overline{H}_8
(\widehat B(\RP^2,4), \pm \Q)$ is trivial because it is Poincar\'e dual to the
0-dimensional homology group with coefficients in a non-constant local system.
Therefore also $\overline{H}_5 (\widehat B(\RP^2,4), \pm \Q) = 0$ by the Euler
characteristic argument. \hfill $\Box$
\medskip

\noindent {\bf Self-join and Caratheodory theorem}. Let us embed a finite cell
complex $M$ generically into the space $\R^T$ of a very large dimension, and
denote by $M^{*r}$ the union of all $(r-1)$-dimensional simplices in $\R^T$,
whose vertices lie in this embedded manifold (and the ``genericity'' of the
embedding means that these simplices have no interior intersections: if two
such simplices have a common point in $\R^T$, then their minimal faces,
containing this point, do coincide).

\begin{proposition}[C.~Caratheodory theorem, see e.g. \cite{hompol}, \cite{KK}]
\label{carat} For any $r \ge 1$, the space $(S^1)^{*r}$ is homeomorphic to
$S^{2r-1}$. \hfill $\Box$ \end{proposition}

\begin{lemma}[see \cite{ser}] \label{lem33}
If $r$ is even, then the sphere $(S^1)^{*r} \sim S^{2r-1}$ has a canonical
orientation, and hence a canonical generator of the $(2r-1)$-dimensional
homology group. If $r$ is odd, then the homeomorphism of $(S^1)^{*r}$ to
itself, induced by any orientation-reversing automorphism of $S^1$, also
violates the orientation. \hfill $\Box$
\end{lemma}

We will use also the following ``homotopy'' modification of Cara\-theo\-dory
theorem.

Let $\mbox{Sym}^i(M)$ be the $i$th symmetric power $M^i/S(i)$ of a manifold
$M$. The disjoint union $\mbox{Sym}^1(M) \sqcup \dots \sqcup \mbox{Sym}^r(M)$
is a partially ordered set by incidence of collections of points. Consider the
join $\mbox{Sym}^1(M) * \dots * \mbox{Sym}^r(M)$ (i.e. the union of all
$(r-1)$-dimensional simplices, whose vertices correspond to the points of
$\mbox{Sym}^1(M), \dots, \mbox{Sym}^r(M)$), and define the topological order
complex of this poset as the union of only those simplices, all whose vertices
are incident to one another. This space can be identified with the space
$\mbox{Sym}^{*r}(M)$, i.e. with the quotient space  $(M * \dots * M)/ S(r)$ of
the join of $r$ copies of $M$ by the obvious action of the symmetric group.
Indeed, for any $j=1, \dots, r$ the space $\mbox{Sym}^j(M)$ can be realised as
the union of barycenters of $(j-1)$-dimensional simplices constituting the join
$M * \dots * M$.

\begin{proposition}[see \cite{KK}]
\label{carat2} For any $M$, there is a canonical homotopy equivalence between
the order complex $\mbox{Sym}^{*r}(M)$ and the space $M^{*r}$ from the previous
lemma. In particular, $\mbox{Sym}^{*r}(S^1)$ is canonically homotopy equivalent
to $S^{2r-1}$. \hfill $\Box$
\end{proposition}

\section{Calculation of the spectral sequence}

\subsection{First column}

The first term $\Psi_1$ of our filtration of $\Psi$ is the space $\RP^2$, so
the column $p=1$ has only one non-trivial cell $E^1_{1,-1} \sim \Q$.

\subsection{Second column}

The term $\Psi_2 \setminus \Psi_1$ consists of two strata, see item $J_2$ in
Lemma \ref{lem3}. The stratum $\Psi_2(b)$ is closed in $\Psi_2 \setminus
\Psi_1$, and we have the exact sequence of the pair
\begin{equation} \label{es2}\dots \to \overline{H}_i(\Psi_2(b), \Q) \to
\overline{H}_i(\Psi_2 \setminus \Psi_1, \Q) \to \overline{H}_i(\Psi_2(a), \Q)
\to \dots
\end{equation}

\begin{lemma}
\label{lem48} Both groups $\overline{H}_*(\Psi_2(b), \Q)$,
$\overline{H}_*(\Psi_2(a), \Q)$ are trivial in all dimensions.
\end{lemma}

{\it Proof.} (a) The stratum $\Psi_2(a)$ is the space of the fiber bundle
(\ref{ltfb}) over the configuration space $B(\RP^2,2)$, i.e. the space of
unordered pairs $\{x \neq y\} \subset \RP^2$. Its fiber $\breve \Psi(\{x \cup
y\})$ is an open interval whose endpoints are associated with the points $x$
and $y$. This interval is constituted by two half-intervals $(\{x\}, \{x \cup
y\}]$ and $(\{y\}, \{x \cup y\}]$ in $\Psi_2$, where $\{x\}$ and $\{y\}$ are
points of $J_1$ corresponding to the points $x, y$, and $\{x \cup y\} \in J_2$.
The endpoints $\{x\}$ and $\{y\}$ are excluded from the segment as they belong
to the first term $\Psi_1$ of our filtration.

Any loop in $B(\RP^2,2)$ violates the orientations of fibers of this bundle if
and only if it permutes these points $x$ and $y$. Thus the $m$-dimensional
rational Borel--Moore homology group of the total space of this fiber bundle is
equal to $\overline{H}_{m-1}(B(\RP^2, 2), \pm \Q)$. By Lemma \ref{lem11} all
such groups are trivial.

(b) The stratum $\Psi_2(b)$ of $\Psi_2 \setminus \Psi_1$ also is the space of a
fiber bundle. Its base is the the projectivized tangent bundle of $\RP^2$, and
the fiber over a tangent element at the point $x \in \RP^2$ is the
half-interval connecting the point of $J_2$ corresponding to this tangent
element with the point $\{x\} \in J_1$. The Borel--Moore homology group of any
fiber bundle with such fibers is trivial in all dimensions. \hfill $\Box$
\medskip

So, both marginal groups in (\ref{es2}) are trivial, hence also the middle
group is, and the column $E^1_{2,q}$ of our spectral sequence is empty.

\subsection{Third column of the spectral sequence}
\label{col3}

The space $\Psi_3 \setminus \Psi_2$ consists of four strata, see Lemma
\ref{lem3}. Let us filter it by
\begin{equation} \label{fltr3} \begin{array}{c}\Psi_3(d) \subset (\Psi_3(d) \cup
\Psi_3(c)) \subset (\Psi_3(d) \cup \Psi_3(c) \cup \Psi_3(b)) \subset \\
\subset (\Psi_3(d) \cup \Psi_3(c) \cup \Psi_3(b) \cup \Psi_3(a))
\end{array}
\end{equation} and study the
auxiliary spectral sequence calculating the Borel--Moore homology of $\Psi_3
\setminus \Psi_2$ and generated by this filtration.

\begin{lemma}
\label{lem3cd} The groups $\overline{H}_*(\Psi_3(d),\Q)$ and
$\overline{H}_*(\Psi_3(c),\Q)$ are trivial in all dimensions.
\end{lemma}

{\it Proof.} In the case of the stratum $\Psi_3(d)$, the link $\partial
\Psi(X)$ of the fiber $\Psi(X)$, $X\in J_3(d),$ is homeomorphic to the closed
disc. Indeed, let $x$ be the point of $\RP^2$ defined by $X$, $\{x\}$ the
corresponding point of $J_1$, and $PT_x(\RP^2) \sim S^1$ the set of tangent
directions at this point, considered as points of the stratum $\Psi_2(b)$ in
$J_2$. Then our link $\partial \Psi(X)$ consists of all segments in $\Psi_2$
connecting the point $\{x\}$ with all points of the projective line
$PT_x(\RP^2)$. Therefore by (\ref{esbh}) the Borel--Moore homology group of the
fiber $\breve \Psi(X)$ is zero in all dimensions, and the same is true for the
entire stratum $\Psi_3(d)$.

For $X \in J_3(c)$, the link $\partial \Psi(X)$ is homeomorphic to a closed
line segment. Indeed, the configuration $X \in J_3$ (consisting of two points
$x \neq y \in \RP^2$ and a tangent direction $l$ at $x$) has two incident
points $\{x\}$ and $\{y\}$ in $J_1$, two incident points $(x,l)$ and $\{x \cup
y\}$ in $J_2$, and three segments in $\partial \Psi(X)$ connecting $\{x\}$ with
$(x,l)$, $\{x\}$ with $\{x \cup y\}$, and $\{y\}$ with $\{x \cup y\}$. This
proves the assertion of Lemma \ref{lem3cd} concerning the homology of
$\Psi_3(c)$. \hfill $\Box$.

\begin{lemma}
\label{lem39} $\overline{H}_i(\Psi_3(b), \Q) \simeq
\overline{H}_{i-4}(\RP^2,\Q)$. In particular, it is equal to $\Q$ if $i=4$ and
is trivial in all other dimensions.
\end{lemma}

{\it Proof.} The stratum $\Psi_3(b)$ is fibered over the space
$\widehat{\RP^2}$ of all lines in $\RP^2$. By Proposition \ref{carat2} the link
$\partial \Psi(L)$ of its fiber over any line $L \in \widehat{\RP^2}$ is
homotopy equivalent to $S^3$, so that the Borel--Moore homology group of
$\breve \Psi(X)$ is equal to $\Q$ in dimension 4 and is trivial in all other
dimensions. By Lemma \ref{lem33}, the monodromy over loops in $\widehat{\RP^2}$
acts trivially on this homology group. Therefore our lemma follows from the
spectral sequence of this fiber bundle. \hfill $\Box$

\begin{lemma}
\label{lem35} $\overline{H}_i(\Psi_3(a), \Q) = \overline{H}_{i-2}(\widehat
B(\RP^2,3), \pm \Q).$ In particular $($by statement 2 of Corollary \ref{lem17}
of Lemma \ref{lem12} and statement 2 of Lemma \ref{lem13}$)$ it is equal to
$\Q$ if $i=8$ or $i=5$ and is trivial in all other dimensions.
\end{lemma}

{\it Proof.} The stratum $\Psi_3(a)$ is fibered (\ref{ltfb}) over the
configuration space $\widehat B(\RP^2,3)$. Its fiber over a configuration is
equal to the open triangle (whose vertices are associated with the points of
the configuration, and the entire triangle is the order complex of non-empty
subsets of this configuration). Any generator of $\pi_1(\widehat B(\RP^2,3))$
described in the proof of Corollary \ref{lem17} (i.e. a rotation of a cube)
changes the orientations of this bundle. Therefore our lemma follows from the
spectral sequence of this fiber bundle. \hfill $\Box$
\medskip

Summing up Lemmas \ref{lem3cd}, \ref{lem39} and \ref{lem35}, we see that the
auxiliary spectral sequence calculating $\overline{H}_*(\Psi_3 \setminus
\Psi_2, \Q)$ has reduced to the exact sequence of two right-hand terms in
(\ref{fltr3}). The $8$-dimensional generator from Lemma \ref{lem35} surely
gives us a summand, which justifies the cell $E^1_{3,5} \simeq \Q$ in Fig.
\ref{mss}. In addition, we have only one suspicious fragment of this exact
sequence:
\begin{equation}
\label{es38}
\begin{array}{r}
0 \to \overline{H}_{5}(\Psi_3 \setminus \Psi_2, \Q) \to
\overline{H}_{5}(\Psi_3(a),\Q)
\stackrel{\partial}{\to} \overline{H}_{4}(\Psi_3(b),\Q) \to \\
\to \overline{H}_{4}(\Psi_3 \setminus \Psi_2, \Q) \to 0,
\end{array}
\end{equation}
where both middle groups are equal to $\Q$ by Lemmas \ref{lem35} and
\ref{lem39}.

\begin{lemma}
\label{lem32} The differential $\partial: \overline{H}_{5}(\Psi_3(a),\Q) \to
\overline{H}_{4}(\Psi_3(b),\Q)$ in $($\ref{es38}$)$ is non-trivial.
\end{lemma}

{\it Proof.} Let us construct a cycle generating the group $
\overline{H}_{5}(\Psi_3(a),\Q)$. Consider the 3-dimensional submanifold
${\mathcal M}^3 \subset \widehat B(\RP^2,3)$, consisting of triples $\{x,y,z\}
\subset \RP^2$ such that $x = (1:0:0)$, $y \neq x$ is a point in the equatorial
line $E \subset \RP^2$ spanned by $x$ and $(0:1:0)$, and $z$ is an arbitrary
point not in this line.

\begin{lemma}
\label{lem321} ${\mathcal M}^3$ is homeomorphic to $\R^3$ and $($being oriented
in any way$)$ generates the group $\overline{H}_3(\widehat B(\RP^2,3), \pm
\Q).$
\end{lemma}

{\it Proof.} The first assertion is obvious. By the second statement of
Corollary \ref{lem17} of Lemma \ref{lem12} the Poincar\'e pairing between the
groups $\overline{H}_3(\widehat B(\RP^2,3), \pm \Q) \equiv
\overline{H}_3(\widehat B(\RP^2,3), {\mathbb O}r)$ and $H_3(\widehat
B(\RP^2,3),\Q)$ is well-defined an non-de\-ge\-ne\-rate. By Lemma
\ref{lem13}.2, these groups are one-dimensional, and the latter of them is
generated by the fundamental class of the manifold of all possible triples of
pairwise orthogonal lines. The cell ${\mathcal M}^3$ has exactly one
transversal intersection with this manifold. \hfill $\Box$
\medskip

Let $[{\mathcal M}^3] \subset \Psi_3 \setminus \Psi_2$ be the preimage of
${\mathcal M}^3$ under the $2$-dimensional fiber bundle $\Psi_3(a) \to \widehat
B(\RP^2, 3)$. By the Thom isomorphism, it generates the group
$\overline{H}_{5}(\Psi_3(a),\Q)$. Its closure in $\Psi_3 \setminus \Psi_2$
contains exactly those points of the stratum $\Psi_3(b)$, which constitute the
fiber $\breve \Psi(E)$ over the point $E \in \widehat{\RP^2}$ (corresponding to
the equatorial line) in the fiber bundle $\Psi_3(b) \to \widehat {\RP^2}$. This
fiber $\breve \Psi(E)$ (supplied with an arbitrary orientation) generates the
group $\overline{H}_{4}(\Psi_3(b),\Q)$, therefore we need only to calculate the
incidence coefficient $[[{\mathcal M}^3]:\breve \Psi(E)]$.

Geometrically, any generic point $T \in \breve \Psi(E)$ is approached by the
points of ${\mathcal M}^3$ four times. Indeed, such a point indicates uniquely
some two points $v \neq w \in E$. Namely, by definition $T$ is some interior
point of a triangle in $\Psi_3$ spanned by a point of $J_1$ (defined by some
point of $\RP^2)$, a point of the main stratum $B(\RP^2, 2)$ of $J_2,$ and the
point $\{E\} \in J_3$. The desired points $v$ and $w$ are exactly the points of
the 2-configuration defining the second vertex of this triangle. Since $T$ is
generic, we will assume that neither $v$ nor $w$ coincides with the fixed point
$x = (1:0:0) \in E$. Let us fix an orientation of the line $E$. We will denote
by $v$ (respectively, $w$) the point of our 2-configuration which is next to
$x$ along (respectively, against) this orientation. When we approach such a
point $T$ from the side of ${\mathcal M}^3$, the point $y$ can approach either
$v$ or $w$ along the line $E$, and $z$ can approach the remaining point $w$ or
$v$ in two different ways, in the same way as the one-dimensional cell of the
standard cell decomposition of $\RP^2$ is approached from the two-dimensional
cell. In total this gives us four possibilities.

\begin{lemma}
\label{lemmma} \label{lem333} All four approaches of $\breve \Psi(E)$ from the
side of $[{\mathcal M}^3]$ have equal $($of the same sign$)$ contributions to
their incidence coefficient.
\end{lemma}

{\it Proof.} To study the incidence coefficients of these subvarieties in
$\Psi_3 \setminus \Psi_2$, let us fix their orientations.

The variety $\breve \Psi(E)$ is 4-dimensional, and is fibered over the
configuration space $B(E \setminus \{x\}, 2)$; its fiber over the configuration
$(v,w)$ consists of two triangles
\begin{equation} \label{oc34} (\{v\}, \{v,w\}, \{E\}) \mbox{ \ and \ } (\{w\},
\{v,w\}, \{E\}).
\end{equation}
Let us orient the first triangle by the indicated sequence of vertices, and the
second triangle by the opposite one, so that their sum has no boundary at the
edge $[\{v, w\}, \{E\}]$. Further, we orient the base by the sequence
$(\partial/\partial v, \partial/ \partial w)$, and orient the whole fiber
bundle by this pair of orientations of the base and the fibres.

The subvariety $[{\mathcal M}^3] \subset \Psi_3 \setminus \Psi_2$ is swept out
by order complexes $\Psi(x,y,z)$ for $(x,y,z) \in {\mathcal M}^3$. Any such
order complex consists of six triangles like $(\{x\}, \{x,y\}, \{x,y,z\})$,
$(\{z\}, \{x,z\}, \{x,y,z\})$, etc. So, we have six families of triangles
parameterized by configurations $(x,y,z) \in {\mathcal M}^3$. Only two of these
six 3-parameter families, namely, the ones swept out by triangles
\begin{equation}
\label{oc35} (\{y\}, \{y,z\}, \{x,y,z\}) \mbox{ \ and \ } (\{z\}, \{y,z\},
\{x,y,z\}), \end{equation} can approach the variety $\breve \Psi(E)$ in its
interior points. Therefore it is enough to consider the 5-dimensional
subvariety $[{\mathcal M}^3]^\circ \subset [{\mathcal M}^3]$ fibered over
${\mathcal M}^3$, with the fiber over $(x,y,z)$ equal to the sum of such two
simplices. Again, we orient the first of them by the indicated order of its
vertices, and the second by the opposite one, so that these orientations are
compatible at the common edges $[\{y,z\}, \{x,y,z\}]$ of these triangles. Also,
we orient the manifold ${\mathcal M}^3$ of configurations $(y,z)$ by the pair
\{the chosen orientation of the line $E \ni y$; the orientation of the upper
hemisphere inducing this orientation on its boundary $E$\}. The orientation of
the entire 5-dimensional cell $[ {\mathcal M}^3]^\circ$ consists of this
orientation of the base and the above-defined orientation of the fibers.

When the triple $(x,y,z) \in {\mathcal M}^3$ degenerates (i.e. $z$ tends to the
equatorial line), the point $(x,y,z) \in J_3$ tends to $\{E\}$, and the sum of
two triangles (\ref{oc35}) (taken with above-described orientations) tends to
$\pm$ the sum of triangles (\ref{oc34}), where $v$ and $w$ are the limit
positions of $y$ and $z$ (or $z$ and $y$).

To prove Lemma \ref{lem333} it remains to compare these limit orientations at
the points of $\breve \Psi(E)$ defined by different approaches from the side of
$[{\mathcal M}^3]^\circ$. Consider a point $(v, w) \in B(E \setminus \{x\}, 2)$
and the sum of oriented triangles (\ref{oc34}) over it.

First we compare two adjacencies $[ {\mathcal M}^3]^\circ \to \breve \Psi(E)$,
when $y = v$ and $z$ tends to $w$ in two possible ways. All ingredients of
these adjacencies (concerning the behavior of fibers and the point $y$) except
for these two ways of tending $z$ to $w$ are exactly the same. Therefore the
contributions of these adjacencies are equal to one another by the same reasons
as for the two adjacencies of 2- and 1-dimensional cells of the standard cell
decomposition of $\RP^2$. In exactly the same way, two different adjacencies
when $y=w$ and $z$ tends to $v$ also give equal contributions.

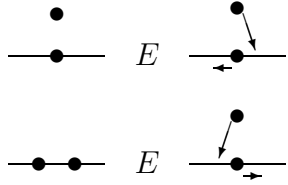
\begin{figure}
\begin{picture}(46,28)
\put(0,2){\line(1,0){16}} \put(5,2){\circle*{2}} \put(11,2){\circle*{2}}
\put(0,20){\line(1,0){16}} \put(8,27){\circle*{2}} \put(8,20){\circle*{2}}
\put(21,18.5){$E$} \put(30,2){\line(1,0){16}} \put(38,2){\circle*{2}}
\put(38,10){\circle*{2}} \put(37,9){\vector(-1,-3){2}}
\put(39,0){\vector(1,0){3}} \put(30,20){\line(1,0){16}}
\put(38,28){\circle*{2}} \put(38,20){\circle*{2}} \put(39,27){\vector(1,-3){2}}
\put(37,18){\vector(-1,0){3}} \put(21,0.5){$E$}
\end{picture}
\caption{Incidence coefficients in the configuration space} \label{conf}
\end{figure}

It remains to compare the contributions of two adjacencies from different
pairs: one with $y$ tending to $v$ and $z$ to $w$, and vice versa for the
other. Some two such adjacencies are shown in Fig. \ref{conf}. The
2-configuration shown in the upper left corner of this picture can degenerate
into the one shown in the lover left corner in two ways shown in the right part
of the picture. Since the orders of two points in $E$ obtained from $y$ and $z$
in these ways are different, the boundary orientations of the base $B(E
\setminus \{x\},2)$ induced by the chosen orientation of ${\mathcal M}^3$ are
opposite. At the same time, the fibers (i.e. the sums of oriented triangles
(\ref{oc35})) in one case approach the corresponding fiber (\ref{oc34})
respecting its chosen orientation, and in the other case reversing this
orientation. Therefore the boundary orientations of the entire fiber bundle
$[E]$ induced from the orientation of $[ {\mathcal M}^3]^\circ$ in these two
cases differ from one another by the factor $(-1) \times (-1)$ (one difference
on the level of bases, and the other by the comparison of the fibers), hence
they also do coincide.

Lemma \ref{lem32} is proved, and column $p=3$ of Fig. \ref{mss} is justified.
\hfill $\Box$

\subsection{Fourth column of the spectral sequence}

The space $\Psi_4 \setminus \Psi_3$ consists of five strata, see Lemma
\ref{lem3}. Let us filter it by
\begin{equation} \begin{array}{l} \label{fltr4} \Psi_4(e) \subset (\Psi_4(e) \cup
\Psi_4(d))
\subset (\Psi_4(e) \cup \Psi_4(d) \cup \Psi_4(c)) \subset \\ \subset (\Psi_4(e)
\cup \Psi_4(d) \cup \Psi_4(c) \cup \Psi_4(b)) \subset \\ \subset (\Psi_4(e)
\cup \Psi_4(d) \cup \Psi_4(c) \cup \Psi_4(b) \cup \Psi_4(a))
\end{array}
\end{equation} and study the auxiliary spectral sequence generated by this
filtration and calculating the group $\overline{H}_*(\Psi_4 \setminus \Psi_3,
\Q)$.

\begin{lemma}
\label{lem442} For any point $X \in J_4$, which belongs to one of strata
described in items $(J_4)(e)$, $(J_4)(d)$ or $(J_4)(c)$ of Lemma \ref{lem3},
the reduced homology group $\tilde{H}_*(\partial \Psi(X),\Q)$ of the link of
the corresponding order complex is trivial in all dimensions. In particular,
the rational Borel--Moore homology groups of corresponding pieces $\Psi_4(e),$
$\Psi_4(d)$ and $\Psi_4(c)$ of $\Psi_4 \setminus \Psi_3$ also are trivial in
all dimensions.
\end{lemma}

{\it Proof.} Any point $X \in J_4(e)$ is defined by a line $L \subset \RP^2$
and a singular point $z \in L$. For any such pair $X=(L,z)$, the link $\partial
\Psi(X)$ is the union of the order complexes $\Psi(L)$ and $\Psi(z)$, where $L$
and $z$ are considered as some points of $J_3$ listed in items $(J_3)(b)$ and
$(J_3)(d)$ of Lemma \ref{lem3}, respectively. The intersection of these two
complexes is the closed segment in $\Psi_2$, connecting the point $\{z\} \in
J_1$ and the point $\{z; \mbox{the direction of }L \mbox{ at } z\}$ of the
stratum of $J_2(b)$. Therefore the assertion of our lemma concerning $(J_4)(e)$
follows from the Mayer-Vietoris sequence of this union.

\unitlength=1.5mm \special{em:linewidth 0.4pt} \linethickness{0.4pt}
\begin{figure}
\mbox{
\begin{picture}(20,15)
\put(5,0){\circle*{1}} \put(15,0){\circle*{1}} \put(4,1){\line(-1,1){3}}
\put(6,1){\line(1,1){3}} \put(14,1){\line(-1,1){3}} \put(16,1){\line(1,1){3}}
\put(0,5){\makebox(0,0)[cc]{$\times$}} \put(9.3,5){\circle*{1}}
\put(10.7,5){\circle*{1}} \put(20,5){\makebox(0,0)[cc]{$\times$}}
\put(1,6){\line(1,1){3}} \put(9,6){\line(-1,1){3}} \put(11,6){\line(1,1){3}}
\put(19,6){\line(-1,1){3}} \put(4.3,10){\makebox(0,0)[cc]{$\times$}}
\put(5.7,10){\circle*{1}} \put(15.7,10){\makebox(0,0)[cc]{$\times$}}
\put(14.3,10){\circle*{1}}
\end{picture}
} \qquad \mbox{
\begin{picture}(30,15)
\put(5,0){\circle*{1}} \put(15,0){\circle*{1}} \put(25,0){\circle*{1}}
\put(4,1){\line(-1,1){3}} \put(7,1){\line(4,1){12}} \put(14,1){\line(-3,1){12}}
\put(14,1){\line(-1,1){3}} \put(16,1){\line(4,1){12}}
\put(24,1){\line(-1,1){3}} \put(26,1){\line(1,1){3}}
\put(0,5){\makebox(0,0)[cc]{$\bullet \bullet$}}
\put(10,5){\makebox(0,0)[cc]{$\times$}} \put(20,5){\makebox(0,0)[cc]{$\bullet \
\bullet$}} \put(30,5){\makebox(0,0)[cc]{$\bullet \bullet$}}
\put(5,10){\makebox(0,0)[cc]{$\bullet \times$}}
\put(15,10){\makebox(0,0)[cc]{$\bullet \bullet \bullet$}}
\put(25,10){\makebox(0,0)[cc]{$\times \bullet$}} \put(1,6){\line(1,1){3}}
\put(1,6){\line(5,1){14}} \put(9,6){\line(-1,1){3}} \put(11,6){\line(4,1){13}}
\put(19,6){\line(-1,1){3}} \put(29,6){\line(-5,1){12}}
\put(29,6){\line(-1,1){3}}
\end{picture}
} \caption{Order complexes for strata $J_4(d)$ and $J_4(c)$} \label{linksdc}
\end{figure}
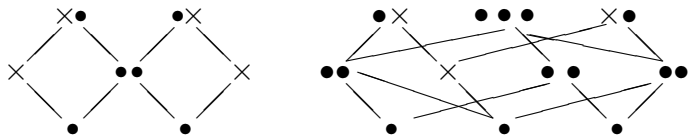

Any point $X \in J_4(d)$ is specified by two different points $ v, w \in \RP^2$
with a generic choice of tangent directions at these points. The link of
$\Psi(X)$ in this case is illustrated in Fig. \ref{linksdc} left. Namely, this
link has two vertices in $J_1$ (defined by these points $v$ and $w$ with
forgotten directions; in the picture they are shown by dots, while the crosses
denote the points supplied with tangent directions), three vertices in $J_2$
(the point $v$ together with the chosen direction, the point $w$ with the
direction, and the configuration $\{v, w\}$ with forgotten directions); and two
points in $J_3$ (configurations $\{v,w\}$ with the direction at only one of
these points). The segments in our picture connect the incident vertices from
$J_1$ and $J_2$ or from $J_2$ and $J_3$. It easy to see that this link is a
contractible simplicial complex, which proves our lemma for the case of
$J_4(d)$.

For $J_4(c)$, the 1-skeleton of the similar link is shown in Fig. \ref{linksdc}
right; the point supplied with a tangent direction corresponds to the middle
point in the bottom row. Again, it is easy to check that the simplicial complex
defined by the poset shown in the picture is contractible. \hfill $\Box$
\medskip

The calculation of the fourth column has reduced to the study of only two
subspaces, $\Psi_4(b)$ and $\Psi_4(a)$; namely, to the exact sequence
\begin{equation} \label{es4} \dots \to \overline{H}_i (\Psi_4(b), \Q) \to
\overline{H}_i (\Psi_4 \setminus \Psi_3, \Q) \to \overline{H}_i (\Psi_4(a), \Q)
\to \overline{H}_{i-1} (\Psi_4(b), \Q) \to \dots
\end{equation}

\begin{lemma} \label{lem4b1}
$\overline{H}_i(\Psi_4(b),\Q) \simeq \overline{H}_{i-7}(\RP^2, {\mathbb O}r)$,
in particular this group is equal to $\Q$ if $i=9$ and is trivial otherwise.
\end{lemma}

{\it Proof.} $\Psi_4(b)$ is fibered over the space $J_4(b)$ of configurations
\begin{equation}
\label{conf4b} \{\mbox{a line } L \in \widehat{\RP^2}; \mbox{a point } z \in
\RP^2 \setminus L\}.\end{equation} The space of such configurations is
obviously diffeomorphic to $T_* \RP^2$, in particular is orientable. Now we
need the following lemma.

\begin{lemma}
\label{lem417} 1. The link $\partial \Psi(X)$ of $\Psi(X)$, $X = (L,z) \in
J_4(b)$, is homotopy equivalent to $S^4$, in particular the group
$\overline{H}_i(\breve \Psi(X), \Q)$ is isomorphic to $\Q$ if $i=5$ and is
trivial for all other $i$.

2. The fundamental group of the space of configurations $($\ref{conf4b}$)$ acts
trivially on the homology bundle associated with the fibers $\breve \Psi(X)$.
\end{lemma}

{\it Proof.} 1. The link $\partial \Psi(X),$ $X=\{L,z\}$, consists of two
subspaces: $A=$ the order complex $\Psi(L),$ and $B=$ the union of all order
complexes $\Psi(\{v,w,z\})$, where $v$ and $w$ are arbitrary (maybe coinciding)
points of the line $L$. The intersection $A \cap B$ is the space $\Psi(L) \cap
\Psi_2$, i.e. the union of similar order subcomplexes $\Psi(\{v,w\})$; by
Proposition \ref{carat2} it is homotopy equivalent to $S^3$. The space $A$ is
compact and contractible, hence the pair $(A, A \cap B)$ is homotopy equivalent
to the pair $(\mbox{cone over } A \cap B, \mbox{its base } A \cap B)$. As for
the space $B$, any of order complexes $\Psi(\{v, w, z\})$ constituting it is
canonically PL-homeomorphic to the triangle with vertices $v, w$ and $z,$ and
hence to the cone over $\Psi(\{v, w\})$ with vertex $\{z\}$. These
homeomorphisms over all pairs $\{v, w\} \in \mbox{Sym}^2(L)$ define a
homeomorphism between the union $B$ of all these triangles and another cone
over $\partial \Psi(L) \cap \Psi_2$. Therefore our space $A \cup B$ is homotopy
equivalent to the suspension over $A \cap B \sim S^3$.

2. By Lemma \ref{lem33}, the generator of $\pi_1(\RP^2)$ acts trivially on the
group $H_3(\Psi(L) \cap \Psi_2, \Q)$. Two parts $A$ and $B$ of $\partial
\Psi(\{L,z\}) \sim \Sigma (\Psi(L) \cap \Psi_2)$ are of different origin and
cannot be permuted by the monodromy, hence the group $H_4(\partial
\Psi(\{L,z\}), \Q)$ also is preserved by this generator. By the natural
boundary isomorphisms (\ref{esbh}), the same is true for the groups
$\overline{H}_*(\breve \Psi(\{L,z\}),\Q)$. \hfill $\Box$
\medskip

Lemma \ref{lem417} is completely proved. Lemma \ref{lem4b1} follows from it
immediately by the composition of the Thom isomorphism of the corresponding
$5$-dimensional fiber bundle (\ref{ltfb}) over the base $J_4(b)$ and the Thom
isomorphism for the non-orientable fibration of this base over
$\widehat{\RP^2}$. \hfill $\Box$
\medskip

\begin{lemma} \label{lem4a1} Any group $\overline{H}_i(\Psi_4(a), \Q)$ is
isomorphic to $\overline{H}_{i-3}(\widehat B(\RP^2,4), \pm \Q)$. In particular,
by Lemma \ref{lem23}, $\overline{H}_*(\Psi_4(a), \Q) \equiv 0$.
\end{lemma}

{\it Proof.} The space $J_4(a)$ is equal to $\widehat B(\RP^2,4)$, and the
fiber of the fiber bundle $\Psi_4(a) \to J_4(a)$ over the configuration $\{x,
y, z, u\} \in \widehat B(\RP^2,4)$  is an open 3-dimensional simplex, whose
vertices correspond to the points of this configuration. Any element of
$\pi_1(\widehat B(\RP^2,4))$ defining an odd permutation of four points
violates the orientation of the bundle of 3-dimensional simplices. Therefore
our lemma follows from the Thom isomorphism for this fiber bundle. \hfill
$\Box$ \medskip

Summing up Lemmas \ref{lem4a1}, \ref{lem4b1} and \ref{lem442}, we obtain that
the auxiliary spectral sequence, calculating the group $\overline{H}_*(\Psi_4
\setminus \Psi_3, \Q)$ and generated by the filtration (\ref{fltr4}), has only
one non-trivial cell, isomorphic to $\Q$ and having total dimension equal to
$9$. This proves the column $p=4$ of Fig. \ref{mss}.

\subsection{The differential $d^1: E^1_{4,5} \to E^1_{3,5}$ is an isomorphism}

Recall that the group $E^1_{4,5}$ is generated by the basic cycle of
$\overline{H}_9(\Psi_4(b), \Q)$, and the group $E^1_{3,5} \equiv
\overline{H}_8(\Psi_3 \setminus \Psi_2, \Q)$ is generated by the fundamental
cycle of the fiber bundle over $\widehat{B}(\RP^2,3)$, whose fiber over the
triple $\{u,v,w\} \subset \RP^2$ is an open triangle with vertices associated
with these points $u, v$ and $w$. Any point of the latter 8-dimensional
manifold participates exactly three times in the links $\partial \Psi(X)$ of
some configurations $X \in J_4(b)$: namely, of the configurations $(L,z)$ where
$z$ is one of our points $u, v$ or $w$, and $L$ is the line spanned by the
remaining two points. Therefore the incidence coefficient of manifolds
generating the groups in question is an odd number, in particular is not equal
to zero. \hfill $\Box$

\subsection{Fifth column of the spectral sequence}

The space $\Psi_5 \setminus \Psi_4$ consists of three strata, see Lemma
\ref{lem3}. Let us filter it by
\begin{equation} \label{fltr5} \Psi_5(c) \subset (\Psi_5(c) \cup \Psi_5(b))
\subset (\Psi_5(c) \cup \Psi_5(b) \cup \Psi_5(a))
\end{equation} and study the auxiliary spectral sequence generated by this
filtration and calculating the Borel--Moore homology of $\Psi_5 \setminus
\Psi_4$.

\begin{lemma} \label{lem5a}
The group $\overline{H}_i(\Psi_5(a), \Q)$ is isomorphic to
$\overline{H}_{i-11}(\RP^2, {\mathbb O}r)$, i.e. is equal to $\Q$ if $i=13$ and
is trivial otherwise.
\end{lemma}

{\it Proof.} The space $J_5(a)$ of non-empty non-singular conics in $\R^3$ is
diffeomorphic to the space of a 3-dimensional vector bundle over $\RP^2$,
namely of the symmetric square of the cotangent bundle: $J_5(a) \simeq S^2
T^*(\RP^2)$, see e.g. Lemma 16 in \cite{ser}. In particular, it is an
orientable 5-dimensional manifold homotopy equivalent to $\RP^2$.

The space $\Psi_5(a) \subset \Psi_5 \setminus \Psi_4$ is fibered over the space
$J_5(a)$ of such conics. For any conic $C \in J_5(a),$ the link $\partial
\Psi(C)$ of the fiber $\breve \Psi(C)$ of this fiber bundle coincides with the
space $\mbox{Sym}^{*4}(S^1)$, see Proposition \ref{carat2}. Therefore by
Propositions \ref{carat} and \ref{carat2} the Borel--Moore homology groups of
the fibers $\breve \Psi(C)$ are isomorphic to $\Q$ in dimension $8$, and are
trivial in all other dimensions. By Lemma \ref{lem33} the fundamental group of
the base acts trivially on these homology groups.

Now, our lemma follows from the Thom isomorphism $\overline{H}_i(\Psi_5(a),\Q)
\simeq \overline{H}_{i-8}(J_5(a), \Q)$ of our orientable fiber bundle over the
space of conics, and the Thom isomorphism $\overline{H}_j(J_5(a), \Q) \simeq
H_{j-3}(\RP^2, {\mathbb O}r)$ for the fibration of this space over $\RP^2$.
\hfill $\Box$
\medskip

Now, let us study $\Psi_5(b)$. The base $J_5(b)$ of the corresponding fiber
bundle (\ref{ltfb}) is the space $B(\widehat{\RP^2},2)$ of pairs of different
lines in $\RP^2$.

\begin{lemma} \label{lem5b}
For any point $X \in B(\widehat{\RP^2},2)$, the link $\partial \Psi(X)$ of the
order complex $\Psi(X)$ is homology equivalent to the sphere $S^7$.
\end{lemma}

{\it Proof.} Let us construct a convenient subdivision of the link $\partial
\Psi(X)$. For any of strata $J_i(\alpha)$, $i <5$, we take all points $Y$ of
this stratum subordinate to $X$ and consider the union of the corresponding
complexes $\breve \Psi(Y)$. This union forms a fiber bundle over the set of
such subordinate points $Y$. Then we calculate Borel--Moore homology groups of
all these fiber bundles; they form the first term of a spectral sequence
calculating the homology of entire $\partial \Psi(X)$ defined by our
stratification.

We can ignore all strata for which the homology group of the fiber $\breve
\Psi(Y)$ is trivial in all dimensions; as was shown previously, it are the
strata $J_2(b),$ $J_3(c),$ $J_3(d),$ $J_4(c),$ $J_4(d),$ and $J_4(e)$. The
remaining configurations subordinated to the point $X$ (i.e. to a pair of
different lines in $\RP^2$) are as follows.
\smallskip

{\bf $J_4(a)$} Two points in one line and two points in the other; none of them
coincides with the intersection point of these lines. The space of such
configurations is equal to $B(\R^1,2) \times B(\R^1,2) \simeq \R^4$, the fiber
$\breve \Psi(Y)$ is a 3-dimensional open simplex. So, the group
$\overline{H}_i$ of this space is equal to $\Q$ for $i=7$ and is trivial in all
other dimensions.
\smallskip

{\bf $J_4(b)$} One of our two lines and one point in the other (different from
their intersection point). The space of such configurations is the union of two
affine lines, and the fiber $\breve \Psi(Y)$ has the Borel--Moore homology
group of $\R^5$, see Lemma \ref{lem417}. So, the Borel--Moore homology group of
entire this space is $\Q^2$ in dimension 6 and is trivial in all other
dimensions.
\smallskip

{\bf $J_3(a)$} The configuration space consists of two three-dimensional cells
and one two-dimensional. Indeed, some two points of our 3-con\-fi\-gu\-ration
should lie in both our lines outside the intersection point, and the third
point can lie also in only one of them or to coincide with the intersection
point. The first two configuration spaces are equal to $B(\R^1,2) \times \R^1
\simeq \R^3$, and the third one to $\R^2$. The fiber $\breve \Psi(Y)$ in all
cases is homeomorphic to an open triangle. Thus, the preimages of the first two
pieces in $\Psi_3(a)$ are homeomorphic to $\R^5$, and that of the third one is
homeomorphic to $\R^4$.
\smallskip

{\bf $J_3(b)$} The configuration space consists of two points (i.e. the choice
of one of two lines constituting $X$). The link $\partial \Psi(Y)$ is homotopy
equivalent to $S^3$, thus the fiber $\breve \Psi(Y)$ has Borel--Moore homology
of $\R^4$.
\smallskip

{\bf $J_2(a)$} Five different strata: both these two points can lie on one or
the other of our two lines (outside the intersection point), or on different
lines, or one of them can coincide with the intersection point, and the other
one lie on one of remaining affine lines. The fiber $\breve \Psi(Y)$ in all
cases is equal to an open interval.
\smallskip

{\bf $J_1$} Three obvious strata equal to $\R^1, \R^1$, and a point; the fiber
in all cases consists of one point.
\smallskip

It is easy to calculate, that the incidence coefficients of two 6-dimen\-sional
cells listed in item $J_4(b)$ and two 5-dimensional cells listed in $J_3(a)$
form the unit matrix. The same is true for the incidence coefficients of two
4-dimensional cells listed in $J_3(b)$ and two first 3-dimensional cells in
$J_2(a)$. The union of remaining cells of dimensions $\le 4$ is a compact
contractible space. Indeed, it is the union of all closed triangles, whose
vertices correspond to the crossing point of our two lines, and two other
points of these lines (one point on each). All these simplices have a common
vertex, hence their union is indeed contractible. So, all cells of our
decomposition of $\partial \Psi(X)$ have killed one another, except only for
the 7-dimensional cell $\Psi_4(a)$ and the 0-dimensional one. \hfill $\Box$

\begin{lemma} \label{lem5b2}
The group $\pi_1(B(\widehat{\RP^2},2))$ acts trivially on the rational homology
groups of the fibers $\partial \Psi(X)$ and $($which is the same by the natural
boundary isomorphisms $($\ref{esbh}$))$ on the rational Borel--Moore homology
groups of fibers $\breve \Psi(X)$.
\end{lemma}

{\it Proof.} Let us define an invariant orientation of the homology bundle with
fiber $H_7(\partial \Psi(X),\Q)$. This group is generated by the fundamental
class of the cell studied in item {\bf $J_4(a)$} of the proof of Lemma
\ref{lem5b}. To orient this cell, we need to choose the orientations of its
base $B(\R^1,2) \times B(\R^1,2)$ and, for any point of this base (i.e. a
configuration of four points), of the fiber over it, i.e. the 3-dimensional
simplex whose four vertices correspond to the points of this configuration. To
do it, we fix an arbitrary numbering of these four points, and define the $j$th
basic tangent vector to the base, $j =1, \dots, 4$, as a shift of the $j$th
point in the direction from the other point of the configuration in the same
line and towards the intersection point of these lines. The orientation of the
fiber (i.e. a tetrahedron) is defined by the same numbering of its vertices. A
different numbering of points simultaneously changes or preserves the
orientations of the base and the fiber, and hence preserves the orientation of
the total space.

The orientation thus defined is obviously invariant and defines a
trivialization of the homology bundle over $B(\widehat{\RP^2},2)$ with the
fiber $H_7(\partial \Psi(X),\Q)$. \hfill $\Box$
\medskip

\begin{corollary} \label{cor5b}
$\overline{H}_i(\Psi_5(b),\Q)$ is equal to
$\overline{H}_{i-8}(B(\widehat{\RP^2},2), \Q)$ for any $i$. In particular, this
group is trivial in all dimensions.
\end{corollary}

{\it Proof.} The first statement follows from Lemmas \ref{lem5b}, \ref{lem5b2}
and the spectral sequence of the corresponding fiber bundle (\ref{ltfb}); the
second statement follows from Lemma \ref{lem11}. \hfill $\Box$
\medskip

\begin{lemma} \label{lem5c}
For any point $X \in J_5(c)$ $($i.e. a  projective line $L$ of multiplicity
2$)$ the homology group of the link $\partial \Psi(X)$ of its order complex is
trivial in all positive dimensions. In particular, the Borel--Moore homology
group of the stratum $\Psi_5(c)$ also is trivial in all dimensions.
\end{lemma}

{\it Proof.} We proceed as in the proof of Lemma \ref{lem5b}. There are only
the following three strata subordinate to $X$ (except for those consisting of
points $Y$ such that $\overline{H}_*(\breve \Psi(Y)) \equiv 0$). In {\bf
$J_3(b)$}, the configuration space of subordinate points consists of only one
point: our line $L \sim \RP^1$ taken once, and the fiber is $\breve \Psi (L)$.
In {\bf $J_2(a)$}, the configuration space is $B(L,2)$, and the fiber is equal
to an 1-dimensional open interval. In {\bf $J_1$}, the configuration space is
our line $L$, the fiber is one point. Also, we are free not to exclude the
stratum $\Psi(X) \cap J_2(b)$, which is fibered over $L$ with the fibers equal
to half-open intervals and anyway does not contribute to the homology group of
$\partial \Psi(X)$.

The union of these four strata is nothing else than the entire order complex
$\Psi(L)$ of the point $L \in J_3$, which is a contractible compact space.
\hfill $\Box$
\medskip

Summing up Lemmas \ref{lem5a}, \ref{lem5c} and Corollary \ref{cor5b}, we obtain
that the auxiliary spectral sequence, calculating the rational Borel--Moore
homology of $\Psi_5 \setminus \Psi_4$ and generated by the filtration
(\ref{fltr5}), has only one non-trivial cell; namely, it is equal to $\Q$ and
its total dimension is equal to $13$. This justifies column $p=5$ of Fig.
\ref{mss}. Theorem \ref{prop21} is completely proved. \hfill $\Box$

\end{document}